\documentclass[11pt]{article}

\usepackage{graphicx,lscape}
\usepackage[mathscr]{euscript}
\usepackage{amsthm,amsmath,amsfonts,amssymb,eepic}
\newcommand{\M}{\mathsf{Mat}}
\newcommand{\B}{\mathcal{B}}
\newcommand{\R}{\mathbb{R}}
\newcommand{\Z}{\mathbb{Z}}

\renewcommand{\int}{\mathop{\rm int}\nolimits}

\newcommand{\conv}{\mathop{\rm convex}\nolimits}

\newtheorem{theorem}{Theorem}[section]
\newtheorem{proposition}[theorem]{Proposition}
\newtheorem{lemma}[theorem]{Lemma}
\newtheorem{corollary}[theorem]{Corollary}
\theoremstyle{definition}
\newtheorem{definition}[theorem]{Definition}
\newtheorem{example}[theorem]{Example}

\theoremstyle{remark}

\numberwithin{equation}{section}

\begin{document}
\title{\textsf{Valuations for matroid polytope subdivisions.}}
\author{\textsf{Federico Ardila\footnote{\textsf{San Francisco State University, San Francisco, CA, USA, federico@math.sfsu.edu.}}} 
\qquad \textsf{Alex Fink\footnote{\textsf{University of California, Berkeley, Berkeley, CA, USA, finka@math.berkeley.edu.}}}
\qquad \textsf{Felipe Rinc\'on\footnote{\textsf{Universidad de Los Andes, Bogot\'a, Colombia,  ed-rinco@uniandes.edu.co.}}}}
\date{}
\maketitle

\begin{abstract}
We prove that the ranks of the subsets and the activities of the bases of a matroid define valuations for the subdivisions of a matroid polytope into smaller matroid polytopes. 
\end{abstract}

\section{\textsf{Introduction.}}

Aside from its wide applicability in many areas of mathematics, one of the pleasant features of matroid theory is the availability of a vast number of equivalent points of view. Among many others, one can think of a matroid as a notion of independence, a closure relation, or a lattice. One point of view has gained prominence due to its applications in algebraic geometry, combinatorial optimization, and Coxeter group theory: that of a matroid as a polytope. This paper is devoted to the study of functions of a matroid which are amenable to this point of view.

To each matroid $M$ one can associate a (basis) \emph{matroid polytope} $Q(M)$, which is the convex hull of the indicator vectors of the bases of $M$. One can recover $M$ from $Q(M)$, and in certain instances $Q(M)$ is the fundamental object that one would like to work with. For instance, matroid polytopes play a crucial role in the matroid stratification of the Grassmannian \cite{GGMS}. They allow us to invoke the machinery of linear programming to study matroid optimization questions \cite{Schrijver}. They are also the key to understanding that matroids are just the type A objects in the family of Coxeter matroids \cite{BGW}.

The subdivisions of a matroid polytope into smaller matroid polytopes have appeared prominently in different contexts: in compactifying the moduli space of hyperplane arrangements (Hacking, Keel and Tevelev \cite{HKT} and Kapranov \cite{Kap}), in compactifying fine Schubert cells in the Grassmannian (Lafforgue \cite{Lafforgue1,Lafforgue2}), and in the study of tropical linear spaces (Speyer \cite{Speyer}). 

Billera, Jia and Reiner \cite{BJR} and Speyer \cite{Speyer, Speyer2}
have shown that some important functions of a matroid, such as its quasisymmetric function and its Tutte polynomial, can be thought of as nice functions of their matroid polytopes: they act as valuations on the subdivisions of a matroid polytope into smaller matroid polytopes. 

The purpose of this paper is to show that two much stronger functions are also valuations.
Consider the matroid functions
\[
f_1(M) = \sum_{A\subseteq[n]} (A, r_M(A)) \,\,\,\,\, \text{and} \,\,\,\,\, 
f_2(M) = \sum_{B \textrm{ basis of } M} (B,E(B),I(B)),
\]
regarded as formal sums. Here $r_M$ denotes matroid rank, and $E(B)$ and $I(B)$ denote the sets of externally and internally active elements of $B$.

\noindent \textbf{Theorems \ref{th:ranks} and \ref{th:activities}}.
\emph{The functions $f_1$ and $f_2$ are valuations for matroid polytope subdivisions: for any subdivision of a matroid polytope $Q(M)$ into smaller matroid polytopes $Q(M_1), \ldots, Q(M_m)$, these functions satisfy
\[ 
 f(M) = \sum_{i} f(M_i) - \sum_{i<j} f(M_{ij}) + \sum_{i<j<k} f(M_{ijk}) - \cdots,
\]
where $M_{ab\ldots c}$ is the matroid whose polytope is $Q(M_a) \cap Q(M_b) \cap \cdots \cap Q(M_c)$.}

The paper is organized as follows.  In Section \ref{sec:prelim} we present some background information on matroids and matroid polytope subdivisions. In Section \ref{sec:valuations} we define valuations under matroid subdivisions, and prove an alternative characterization of them. In Section \ref{sec:family} we present a useful family of valuations, which we use to prove Theorems \ref{th:ranks} and \ref{th:activities} in Section \ref{sec:theorems}.  Finally in Section \ref{sec:related work} we
discuss related work.

\section{\textsf{Preliminaries on matroids and matroid subdivisions.}}\label{sec:prelim}

A \emph{matroid} is a combinatorial object which unifies several notions of independence.  We start with some basic definitions; for more information on matroid theory we refer the reader to \cite{Oxley}.

There are many equivalent ways of defining a matroid. We will adopt the basis point of view, which is the most convenient for the study of matroid polytopes.

\begin{definition}
A \emph{matroid} $M$ is a pair $(E, \B)$ consisting of a finite set $E$ and a collection of subsets $\B$ of $E$, called the \emph{bases} of $M$, which satisfies the \emph{basis exchange axiom:} If $B_1, B_2 \in \B$ and $b_1 \in B_1- B_2$, then there exists $b_2 \in B_2 - B_1$ such that $B_1 - b_1 \cup b_2 \in \B$.
\end{definition}

We will find it convenient to allow $(E, \emptyset)$ to be a matroid; this is not customary.

A subset $A \subseteq E$ is \emph{independent} if it is a subset of a basis. All the maximal independent sets contained in a given set $A \subseteq E$ have the same size, which is called the \emph{rank} $r_M(A)$ of $A$. In particular, all the bases have the same size, which is called the rank $r(M)$ of $M$.

\begin{example}
If $E$ is a finite set of vectors in a vector space, then the maximal linearly independent subsets of $E$ are the bases of a matroid. The matroids arising in this way are called \emph{representable}, and motivate much of the theory of matroids. 
\end{example}

\begin{example}
If $k \leq n$ are positive integers, then the subsets of size $k$ of $[n]=\{1,\ldots, n\}$ are the bases of a matroid, called the \emph{uniform matroid} $U_{k,n}$.
\end{example}

\begin{example}
Given positive integers $1 \leq s_1 < \ldots < s_r \leq n$, the sets $\{a_1, \ldots, a_r\}$ such that $a_1 \leq s_1, \ldots, a_r \leq s_r$ are the bases of a matroid, called the \emph{Schubert matroid} $SM_n(s_1, \ldots, s_r)$. These matroids were discovered by Crapo \cite{C} and rediscovered in various contexts; they have been called shifted matroids \cite{A, K}, PI-matroids \cite{BJR}, generalized Catalan matroids \cite{BDN}, and freedom matroids \cite{CS}, among others. We prefer the name Schubert matroid, which highlights their relationship with the stratification of the Grassmannian into Schubert cells \cite[Section 2.4]{OMbook}.
\end{example}

The following geometric representation of a matroid is central to our study.

\begin{definition} 
Given a matroid $M=([n],\B)$, the (basis) \emph{matroid polytope} $Q(M)$ of $M$ is the convex hull of the indicator vectors of the bases of $M$:
\[
Q(M) = \conv\{e_B \, : \, B \in \B\}.
\]
For any $B=\{b_1,\ldots,b_r\}\subseteq [n]$, 
by $e_B$ we mean $e_{b_1} + \cdots + e_{b_r}$,
where $\{e_1, \ldots, e_n\}$ is the standard basis of $\R^n$.
\end{definition}

When we speak of ``a matroid polytope", we will refer to the polytope of a specific matroid, in its specific position in $\R^n$. The following elegant characterization is due to Gelfand and Serganova:

\begin{theorem}\label{tm:GS} \cite{GGMS} Let $\B$ be a collection of $r$-subsets of $[n]$ and let 
$Q(\B)=\conv\{e_B \, : \, B \in \B\}$. The following are equivalent:

\noindent 1. $\B$ is the collection of bases of a matroid. 

\noindent 2. Every edge of  $Q(\B)$ is a parallel translate of  $e_i-e_j$ for some $i,j \in [n]$.
\end{theorem}

When the statements of Theorem \ref{tm:GS} are satisfied, the edges of $Q(\B)$ correspond exactly to the pairs of different bases $B, B'$ such that $B' = B \setminus i \cup j$ for some $i,j \in [n]$. Two such bases are called {\em adjacent bases}.

\medskip

A {\em subdivision} of a polytope $P$ is a set of polytopes $S=\{P_1, \ldots, P_m\}$, whose vertices are vertices of $P$, such that
\begin{itemize}
\item $P_1 \cup \cdots \cup P_m = P$, and
\item for all $1 \leq i<j \leq m$, if the intersection $P_i \cap P_j$ is nonempty, then it is a proper face of both $P_i$ and $P_j$.
\end{itemize}
The {\em faces} of the subdivision $S$ are the faces of the $P_i$; it is easy to see that the interior faces of $S$ (i.e. faces not contained in the boundary of $P$) are exactly the non-empty intersections between some of the $P_i$.

\begin{definition}
A {\em matroid polytope subdivision} is a subdivision of a matroid polytope $Q=Q(M)$ into matroid polytopes $Q_1=Q(M_1), \ldots, Q_m=Q(M_m)$.  We will also refer to this as a \emph{matroid subdivision} of the matroid $M$ into $M_1, \ldots, M_m$.
\end{definition}

The lower-dimensional faces of the subdivision, which are intersections of subcollections of the $Q_i$, are also of interest.  
Given a set of indices $A=\{a_1,\ldots,a_s\} \subseteq [m]$,
we will write $Q_A=Q_{a_1\cdots a_s}:=\bigcap_{a\in A}Q_a$.
By convention, $Q_\emptyset = Q$.  Since any face of a matroid polytope is itself a matroid polytope, it follows that any nonempty $Q_A$ is the matroid polytope of a matroid, which we denote $M_A$. 

\medskip

\begin{figure}[h]
\centering
\includegraphics[height=7cm]{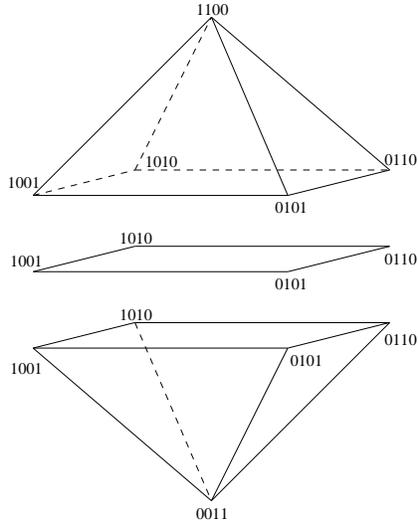}
\caption{The matroid subdivision of a regular octahedron into two square pyramids.}
\label{fig:subdiv}
\end{figure}

Because of the small number of matroid polytopes in low dimensions, there is a general lack of small examples of matroid subdivisions. In two dimensions the only matroid polytopes are the equilateral triangle and the square, which have no nontrivial matroid subdivisions. In three dimensions, the only nontrivial example is the subdivision of a regular octahedron (with bases $\{12,13,14,23,24,34\}$) into two square pyramids (with bases $\{12,13,14,23,24\}$ and $\{13,14,23,24,34\}$, respectively); this subdivision is shown in Figure \ref{fig:subdiv}.

\begin{example}\label{ex:subdiv}
A more interesting example is the following subdivision \cite[Example 7.13]{BJR}:
Let $M_1 = SM_6(2,4,6)$ be the Schubert matroid whose bases are the sets $\{a,b,c\} \subseteq [6]$ such that $a \leq 2, b \leq 4,$ and $c \leq 6$. The permutation $\sigma=345612$ acts on the ground set~$[6]$ of $M_1$, thus defining the matroids $M_2 = \sigma M_1$ and $M_3 = \sigma^2 M_1.$ (Note that $\sigma^3$ is the identity.) Then $\{M_1, M_2, M_3\}$ is a subdivision of $M=U_{3,6}$\footnote{One can easily generalize this construction to obtain a subdivision of $U_{a,ab}$ into $a$ isomorphic matroids.}. 
\end{example}

\section{\textsf{Valuations under matroid subdivisions.}}\label{sec:valuations}

We now turn to the study of matroid functions which are valuations under the subdivisions of a matroid polytope into smaller matroid polytopes. Throughout this section, $\M = \M_n$ will denote the set of matroids with ground set $[n]$, and $G$ will denote an arbitrary abelian group. As before, given a subdivision $\{M_1, \ldots, M_m\}$ of a matroid $M$ and a subset $A \subseteq [m]$, $M_A$ is the matroid whose polytope is $\bigcap_{a \in A} Q(M_a)$.

\begin{definition}
A function $f:\M \rightarrow G$ is a {\em valuation under matroid subdivision}, or simply a {\em valuation}\footnote{This use of the term \emph{valuation} is standard in convex geometry \cite{Mc}. It should not be confused with the unrelated notion of a matroid valuation found in the theory of valuated matroids \cite{DW}.}, if for any subdivision $\{M_1, \ldots, M_m\}$ of a matroid $M \in \M$, we have
\begin{equation}\label{wellbeh}
\sum_{A\subseteq [m]} (-1)^{|A|} f(M_A) = 0
\end{equation}
or, equivalently,
\begin{equation}\label{wellbeh2}
f(M) = \sum_{i} f(M_i) - \sum_{i<j} f(M_{ij}) + \sum_{i<j<k} f(M_{ijk}) - \cdots
\end{equation}
\end{definition}

Recall that, contrary to the usual convention, we have allowed $\emptyset = ([n], \emptyset)$ to be a matroid. We will also adopt the convention that $f(\emptyset)=0$ for all the matroid functions considered in this paper. 

Many important matroid functions are well-behaved under subdivision. Let us start with some easy examples.

\begin{example}
The function $\textsf{Vol}$, which assigns to each matroid $M \in \M$ the $n$-dimensional volume of its polytope $Q(M)$, is a valuation. This is clear since the lower-dimensional faces of a matroid subdivision have volume $0$.
\end{example}

\begin{example}
The \emph{Ehrhart polynomial} $E_P(x)$ of a lattice polytope $P$ in $\R^d$ is the polynomial such that, for a positive integer $n$, $E_P(n)=|nP \cap \Z^d|$ is the number of lattice points contained in the $n$-th dilate $nP$ of $P$ \cite[Section 4.6]{EC1}. By the inclusion-exclusion formula, the function $E:\M \rightarrow \R[x]$ defined by $E(M) = E_{Q(M)}(x)$ is a valuation.
\end{example}

\begin{example}
The function $b(M) = (\textrm{number of bases of }M)$ is a valuation. This follows from the fact that the only lattice points in $Q(M)$ are its vertices, which are the indicator vectors of the bases of $M$; so $b(M)$ is the evaluation of $E(M)$ at $x=1$.
\end{example}

Before encountering other important valuations, let us present an alternative way of characterizing them. 

\begin{theorem}\label{th:equiv}
A function $f: \M \to G$ is a valuation if and only if, for any matroid subdivision $S$ of $Q=Q(M)$, 
\begin{equation}\label{igual}
f(M) = \sum_{F \in \int(S)} (-1)^{\dim(Q)-\dim(F)}f(M(F)),
\end{equation}
where the sum is over the interior faces of the subdivision $S$, and $M(F)$ denotes the matroid whose matroid polytope is $F$.
\end{theorem}

To prove Theorem \ref{th:equiv} we first need to recall some facts from topological combinatorics. These can be found, for instance, in \cite[Section 3.8]{EC1}.

\begin{definition}\label{def:rcc}
A \emph{regular cell complex} is a finite set $C=\left\{ \sigma_1, \sigma_2, \dotsc, \sigma_s \right\}$ of pairwise disjoint and nonempty \emph{cells} $\sigma_i \subseteq \mathbb{R}^d$ such that for any $i \in [s]$:
\begin{enumerate}
\item 
$\overline{\sigma_i} \approx \mathbb{B}^{m_i}$ and $\overline{\sigma_i} \setminus \sigma_i \approx \mathbb{S}^{m_i-1}$ for some nonnegative integer $m_i
$, called the \emph{dimension} of $\sigma_i$.
\item $\overline{\sigma_i} \setminus \sigma_i$ is the union of some other $\sigma_j$s.
\end{enumerate}
Here $\overline{\sigma_i}$ denotes the topological closure of $\sigma_i$ and $\approx$ denotes homeomorphism. Also $\mathbb{B}^l$ and $\mathbb{S}^l$ are the $l$-dimensional closed unit ball and unit sphere, respectively. The \emph{underlying space} $\left|C\right|$ of $C$ is the topological space $\sigma_1 \cup \cdots \cup \sigma_s$.
\end{definition}


\begin{definition}
Let $C$ be a regular cell complex, and let $c_i$ be the number of $i$-dimensional cells of $C$. The \textit{Euler characteristic} of $C$ is: 
\[
\chi(C)=\sum_{\sigma \in C} (-1)^{\dim(\sigma)} = \sum_{i \in \mathbb{N}} (-1)^i c_i = c_0-c_1+c_2-c_3 \dotsb.
\]
The \emph{reduced Euler characteristic} of $C$ is  $\tilde{\chi}(C)=\chi(C)-1$. A fundamental fact from algebraic topology is that the Euler characteristic of $C$ depends solely on the homotopy type of the underlying space $\left|C\right|$.
\end{definition}

\begin{definition}
For a regular cell complex $C$, let $P(C)$ be the poset of cells of $C$, ordered by $\sigma_i \leq \sigma_j$ if $\overline{\sigma_i} \subseteq \overline{\sigma_j}$. Let $\widehat{P}(C)=P(C) \cup \{\hat{0},\hat{1}\}$ be obtained from $P(C)$ by adding a minimum and a maximum element. 
\end{definition}

\begin{definition}
The \emph{M\"obius function} $\mu:\textrm{Int}(P) \rightarrow \Z$ of a poset $P$ assigns an integer to each closed interval of $P$, defined recursively by
\[
\mu_P(x,x) = 1, \qquad \sum_{x \leq a \leq y} \mu(x,a) = 0 \,\,\, \textrm{for all }x<y.
\]
It can equivalently be defined in the following dual way:
\[
\mu_P(x,x) = 1, \qquad \sum_{x \leq a \leq y} \mu(a,y) = 0 \,\,\, \textrm{for all }x<y.
\]
\end{definition}

The following special case of Rota's Crosscut Theorem is a powerful tool for computing the M\"obius function of a lattice.

\begin{theorem}\cite{Rota}\label{th:crosscut}
Let $L$ be any finite lattice. Then for all $x \in L$,
\[
\mu(\hat{0},x) = \sum_B (-1)^{|B|},
\]
where the sum is over all sets $B$ of atoms of $L$  such that $\bigvee B = x$.
\end{theorem}

Finally, we recall an important theorem which relates the topology and combinatorics of a regular cell complex.
\begin{theorem}\cite[Proposition 3.8.9]{EC1}\label{posettop}
Let $C$ be a regular cell complex such that $\left|C\right|$ is a manifold, with or without boundary. Let $P = \widehat{P}(C)$. Then
\[
\mu_P(x,y)=
\begin{cases}
\tilde{\chi}(\left|C\right|) & \text{if } x=\hat{0} \text{ and } y=\hat{1}, \\
0 & \text{if } x \neq \hat{0}, y= \hat{1}, \text{ and $x$ is on the boundary of } \left|C\right|,  \\
(-1)^{l(x,y)} & \text{otherwise},
\end{cases}
\]
where $l(x,y)$ is the number of elements in a maximal chain from $x$ to $y$.
\end{theorem}

We are now in a position to prove Theorem \ref{th:equiv}.

\begin{proof}[Proof of Theorem \ref{th:equiv}] 
Let $S=\{M_1, \ldots, M_m\}$ be a matroid subdivision of $M$. Let $\{Q_1, \ldots, Q_m\}$ and $Q$ be the corresponding polytopes. Notice that the (relative interiors of the) faces of the subdivision $S$ form a regular cell complex whose underlying space is $Q$. Additionally, the poset $\widehat{P}(S)$ is a lattice, since it has a meet operation $\sigma_i \wedge \sigma_j = \int(\overline{\sigma}_i \cap \overline{\sigma}_j)$ and a maximum element.

We will show that 
\begin{equation}\label{wellbeh3}
\sum_{F \in \int(S)} (-1)^{\dim(Q)-\dim(F)}f(M(F)) =
\sum_{i} f(M_i) - \sum_{i<j} f(M_{ij}) + \sum_{i<j<k} f(M_{ijk}) - \cdots
\end{equation}
which will establish the desired result in view of (\ref{wellbeh2}). In the right hand side, each term is of the form $f(M(F))$ for an interior face $F$ of the subdivision $S$ and moreover, all interior faces $F$ appear. The term $f(M(F))$ appears with coefficient
\[
\sum_{A \subseteq [m] \, : \, M_A = M(F)} (-1)^{|A|+1}.
\]
This is equivalent to summing over the sets of coatoms of the lattice $\widehat{P}(S)$ whose meet is $F$. By Rota's Crosscut Theorem \ref{th:crosscut}, when applied to the poset $\widehat{P}(S)$ turned upside down, this sum equals $-\mu_{\widehat{P}(S)}(F,\hat{1})$. Theorem \ref{posettop} tells us that this is equal to $(-1)^{l(F,\hat{1})-1} = (-1)^{\dim(Q) - \dim(F)}$, as desired.
\end{proof}

\section{\textsf{A powerful family of valuations.}}\label{sec:family}

\begin{definition}
Given $X \subseteq \mathbb{R}^n$, let $i_X: \M \to \mathbb{Z}$ be defined by
\[
i_X(M)=
\begin{cases}
1 & \text{if } Q(M) \cap X \neq \emptyset, \\
0 & \text{otherwise}.
\end{cases}
\]
\end{definition}

Our interest in these functions is that, under certain hypotheses, they are valuations under matroid subdivisions.
They are a powerful family for our purposes because many valuations of interest, in particular those
of Section~\ref{sec:theorems}, can be obtained as 
linear combinations of evaluations of these valuations.

\begin{theorem}\label{convopen}
If $X \subseteq \mathbb{R}^n$ is convex and open, then $i_X$ is a valuation. 
\end{theorem}
\begin{proof}
Let $M \in \M$ be a matroid and $S$ be a subdivision of $Q= Q(M)$. We can assume that $Q \cap X \neq \emptyset$, or else the result is trivial. We can also assume that $X$ is bounded since $i_X = i_{X \cap [0,1]^n}$.

We will first reduce the proof to the case when $X$ is an open polytope in $\R^n$. By the Hahn-Banach separation theorem \cite[Theorem~3.4]{Rudin}, for each face $F$ of $S$ such that $F \cap X = \emptyset$ there exists an open halfspace $H_F$ containing $X$ and disjoint from $F$. Let 
\[
X'= \bigcap_{F \cap X = \emptyset} H_F
\]
be the intersection of these halfspaces. Then $X' \supseteq X$ and $X' \cap F= \emptyset$ for each face $F$ not intersecting $X$, so $i_{X'}$ and $i_X$ agree on all the matroids of this subdivision. If we define $X''$ as the intersection of $X'$ with some open cube containing $Q$, $i_{X''}$ and $i_X$ agree on this subdivision and $X''$ is an open polytope. 

We can therefore assume that $X$ is an open polytope in $\R^n$; in particular it is full-dimensional. Note that $X \cap \int(Q)$ is the interior $\int(R)$ of some polytope $R\subseteq Q$. Since $R$ and $Q$ have the same dimension, $R \approx \mathbb{B}^{\dim(Q)}$ and $\partial R \approx \mathbb{S}^{\dim(Q) - 1}$. If $F$ is a face of the subdivision $S$ and $\sigma$ is a face of the polytope $R$, let $c_{F,\sigma} =  \int(F)\cap\int(\sigma)$. Since $c _{F,\sigma}$ is the interior of a polytope, it is homeomorphic to a closed ball and its boundary to the corresponding sphere. Define
\begin{align*}
C &= \left\{ c_{F,\sigma} : c_{F,\sigma} \neq \emptyset \right\} \\
\partial C &= \left\{ c_{F,\sigma} : c_{F,\sigma} \neq \emptyset \text{ and } \sigma \neq R \right\}.
\end{align*}
The elements of $C$ form a partition of $R$ and in this way $C$ is a regular cell complex whose underlying space is $R$. Similarly, $\partial C$ is a regular subcomplex whose underlying space is $\partial R$. Note that if $F$ is an interior face of $S$, $c _{F,R} = \int(F) \cap \int(R) \neq \emptyset$ if and only if $F \cap X \neq \emptyset$, and in this case $\dim(c _{F,R})=\dim(F)$.

We then have 
\begin{eqnarray*}
 \sum_{F \in \int(S)}(-1)^{\dim(F)}i_X(M(F)) 
&=&  \sum_{\substack{F \in \int(S) \\ F \cap X \neq \emptyset}} 
  (-1)^{\dim(F)} \\
&=& \sum_{c_{F,R}\neq\emptyset} 
  (-1)^{\dim(c_{F,R})} \\
&=& \sum_{c \in C} (-1)^{\dim(c)} - 
  \sum_{c \in \partial C} (-1)^{\dim(c)} \\
&=&  \chi(R) - \chi(\partial R)  \\
&=&  1 - \left( 1 + (-1)^{\dim(Q)-1} \right) \\
&=& (-1)^{\dim(Q)} i_X(M),
\end{eqnarray*}
which finishes the proof in view of Theorem \ref{th:equiv}.
\end{proof}

\begin{corollary}
If $X \subseteq \mathbb{R}^n$ is convex and closed, then $i_X$ is a valuation. 
\end{corollary}
\begin{proof}
As before, we can assume that $X$ is bounded since $i_X = i_{X \cap [0,1]^n}$. Now let $S$ be a subdivision of  $Q=Q(M)$ into $m$ parts. For all $A \subseteq \left[m\right]$ such that $X \cap Q_A = \emptyset$, the distance $d(X, Q_A)$ is positive since $X$ is compact and $Q_A$ is closed. Let  $\epsilon>0$ be smaller than all those distances, and define the convex open set
\[
U = \left\{ x \in \mathbb{R}^n : d(x,X) < \epsilon \right\}.
\]
For all $A \subseteq \left[m\right]$ we have that $X \cap Q_A \neq \emptyset$ if and only if $U \cap Q_A \neq \emptyset$. By Theorem \ref{convopen},
\[
\sum_{A \subseteq \left[m\right]} (-1)^{\left|A \right|}i_X(M_A) = 
\sum_{A \subseteq \left[m\right]} (-1)^{\left|A \right|}i_U(M_A)=0
\]
as desired.
\end{proof}

In particular, $i_P$ is a valuation for any polytope $P \subseteq \mathbb{R}^n$.

\begin{proposition}
The constant function $c(M) = 1$ for $M \in \M$ is a valuation.
\end{proposition}
\begin{proof}
This follows from $c(M) = i_{[0,1]^n}$.
\end{proof}

\begin{proposition}
If $X \subseteq \mathbb{R}^n$ is convex, and is either open or closed, then the function
$\overline{i_X}: \M \to \mathbb{Z}$ defined by
\[
\overline{i_X}(M)=
\begin{cases}
0 & \text{if } Q(M) \cap X \neq \emptyset, \\
1 & \text{otherwise},
\end{cases}
\]
is a valuation.
\end{proposition}
\begin{proof}
Notice that $\overline{i_X}= 1 - i_X$, which is the sum of two valuations.
\end{proof}

\section{\textsf{Subset ranks and basis activities are valuations.}}\label{sec:theorems}

We now show that there are two surprisingly fine valuations of a matroid: the ranks of the subsets and the activities of the bases.

\subsection{\textsf{Rank functions.}}

\begin{theorem}\label{th:ranks}
Let $G$ be the free abelian group on symbols of the form $(A,s)$,
$A\subseteq[n]$, $s\in\mathbb Z_{\geq0}$.  The function
$F:\M \to G$ defined by 
\[
F(M) = \sum_{A\subseteq[n]} (A, r_M(A))
\]
is a valuation.
\end{theorem}

\begin{proof}
It is equivalent to show that the function $f_{A,s}: \M \to \mathbb{Z}$ defined by
\[
f_{A,s}(M)=
\begin{cases}
1 & \text{if } r_{M}(A)=s, \\
0 & \text{otherwise},
\end{cases}
\]
is a valuation. Define the polytope
\[
P_{A,s} = \left\{ x \in \left[0,1\right]^n : \sum_{i \in A} x_i \geq s \right\}. 
\]
A matroid $M$ satisfies that $r_{M}(A)=s$ if and only if it has a basis $B$ with $\left|A \cap B \right| \geq s$, and it has no basis $B$ such that $\left|A \cap B \right| \geq s+1$. This is equivalent to $Q(M) \cap P_{A,s} \neq \emptyset$ and $Q(M) \cap P_{A,s+1} = \emptyset$. It follows that $f_{A,s} = i_{P_{A,s}} - i_{P_{A,s+1}}$, which is the sum of two valuations.
\end{proof}

\subsection{\textsf{Basis activities.}}

One of the most powerful invariants of a matroid is its \emph{Tutte polynomial}:
\[
T_M(x,y) = \sum_{A \subseteq [n]} (x-1)^{r(M)-r(A)} (y-1)^{|A|-r(A)}.
\]
Its importance stems from the fact that many interesting invariants of a matroid satisfy the \emph{deletion-contraction recursion}, and every such invariant is an evaluation of the Tutte polynomial~\cite{Tutte}.

\begin{definition}
Let $B$ be a basis of the matroid $M=([n],\B)$. An element $i \in B$ is said to be \emph{internally active} with respect to $B$ if $i<j$ for all $j \notin B$ such that $B - i \cup j \in \B$. Similarly, an element $i \notin B$ is said to be \emph{externally active} with respect to $B$ if $i<j$ for all $j \in B$ such that $B - j \cup i \in \B$. Let $I(B)$ and $E(B)$ be the sets of internally and externally active elements with respect to $B$.
\end{definition}

\begin{theorem}(Tutte, Crapo \cite{Tutte})\label{tutteacti}
The Tutte polynomial of a matroid is
\[
T_M(x,y) = \sum_{B \textrm{ basis of } M} x^{|I(B)|}y^{|E(B)|}.
\]
\end{theorem}

\begin{theorem} \label{th:activities}
Let $G$ be the free abelian group generated by the triples $(B,E,I)$, where $B \subseteq \left[n\right]$, $E \subseteq \left[n\right] \setminus B$ and $I \subseteq B$. 
The function $F:\M \to G$  defined by
\begin{equation}\label{activity}
F(M) = \sum_{B \textrm{ basis of } M} (B,E(B),I(B))
\end{equation}
is a valuation.
\end{theorem}

Before proving this result, let us illustrate its strength with an example. Consider the subdivision of $M=U_{3,6}$ into three matroids $M_1, M_2,$ and $M_3$ described in Example \ref{ex:subdiv}. Table \ref{tableacti} shows the external and internal activity with respect to each basis in each one of the eight matroids $M_A$ arising in the subdivision. The combinatorics prescribed by Theorem \ref{th:activities} is extremely restrictive: in any row, any choice of $(E,I)$ must appear the same number of times in the $M_A$s with $|A|$ even and in the $M_A$s with $|A|$ odd.

\begin{landscape}
\renewcommand{\tabcolsep}{4pt}
\begin{table}
{ \small
\begin{tabular}{|c||c|c||c|c||c|c||c|c||c|c||c|c||c|c||c|c|} \hline
& \multicolumn{2}{c||}{$M$} & \multicolumn{2}{c||}{$M_1$} & \multicolumn{2}{c||}{$M_2$} & \multicolumn{2}{c||}{$M_{1,2}$} & \multicolumn{2}{c||}{$M_3$} & \multicolumn{2}{c||}{$M_{1,3}$} & \multicolumn{2}{c||}{$M_{2,3}$} & \multicolumn{2}{c|}{$M_{1,2,3}$} \\ \hline
$B$ & $E(B)$ & $I(B)$ & $E(B)$ & $I(B)$ & $E(B)$ & $I(B)$ & $E(B)$ & $I(B)$ & $E(B)$ & $I(B)$ & $E(B)$ & $I(B)$ & $E(B)$ & $I(B)$ & $E(B)$ & $I(B)$ \\ \hline \hline
123 & $\emptyset$ & 123	& $\emptyset$ & 123 & &				& & 			& & 			& & 			& & 			& & 			\\ \hline
124 & $\emptyset$ & 12 	& $\emptyset$ & 12 	& & 			& & 			& & 			& & 			& & 			& & 			\\ \hline
125 & $\emptyset$ & 12 	& $\emptyset$ & 12 	& & 			& & 			& $\emptyset$ & 125 & $\emptyset$ & 125 & & 			& & 			\\ \hline
126 & $\emptyset$ & 12 	& 5 & 12 	& & 			& & 			& $\emptyset$ & 12 	& 5 & 12  & & 			& & 			\\ \hline
134 & $\emptyset$ & 1 	& $\emptyset$ & 1 	& $\emptyset$ & 134	& $\emptyset$ & 134 & & 			& & 			& & 			& & 			\\ \hline
135 & $\emptyset$ & 1 	& $\emptyset$ & 1 	& $\emptyset$ & 13 	& $\emptyset$ & 13 	& $\emptyset$ & 15 	& $\emptyset$ & 15 	& $\emptyset$ & 135 & $\emptyset$ & 135 \\ \hline
136 & $\emptyset$ & 1 	& 5 & 1 	& $\emptyset$ & 13 	& 5 & 13 	& $\emptyset$ & 1 	& 5 & 1 	& $\emptyset$ & 13 	& 5 & 13 	\\ \hline
145 & $\emptyset$ & 1 	& 3 & 1 	& $\emptyset$ & 1 	& 3 & 1 	& 3 & 15 	& 3 & 15 	& 3 & 15 	& 3 & 15 	\\ \hline
146 & $\emptyset$ & 1 	& 35 & 1 	& $\emptyset$ & 1 	& 35 & 1 	& 3 & 1 	& 35 & 1 	& 3 & 1 	& 35 & 1 	\\ \hline
156 & $\emptyset$ & 1 	& & 			& & 			& & 			& $\emptyset$ & 1 	& & 			& & 			& & 			\\ \hline
234 & 1 & $\emptyset$ 	& 1 & $\emptyset$ 	& 1 & 34 	& 1 & 34 	& & 			& & 			& & 			& & 			\\ \hline
235 & 1 & $\emptyset$ 	& 1 & $\emptyset$		& 1 & 3 	& 1 & 3 	& 1 & 5 	& 1 & 5 	& 1 & 35 	& 1 & 35 	\\ \hline
236 & 1 & $\emptyset$ 	& 15 & $\emptyset$ 	& 1 & 3 	& 15 & 3 	& 1 & $\emptyset$ 	& 15 & $\emptyset$ 	& 1 & 3 	& 15 & 3 	\\ \hline
245 & 1 & $\emptyset$ 	& 13 & $\emptyset$ 	& 1 & $\emptyset$ 	& 13 & $\emptyset$ 	& 13 & 5 	& 13 & 5 	& 13 & 5 	& 13 & 5 	\\ \hline
246 & 1 & $\emptyset$ 	& 135 & $\emptyset$ & 1 & $\emptyset$ 	& 135 & $\emptyset$ & 13 & $\emptyset$ 	& 135 & $\emptyset$ & 13 & $\emptyset$ 	& 135 & $\emptyset$ \\ \hline
256 & 1 & $\emptyset$ 	& & 			& & 			& & 			& 1 & $\emptyset$ 	& & 			& & 			& & 			\\ \hline
345 & 12 & $\emptyset$	& & 			& 12 & $\emptyset$ 	& & 			& & 			& & 			& & 			& & 			\\ \hline
346 & 12 & $\emptyset$	& & 			& 12 & $\emptyset$ 	& & 			& & 			& & 			& & 			& & 			\\ \hline
356 & 12 & $\emptyset$	& & 			& 12 & 3 	& & 			& 12 & $\emptyset$ 	& & 			& 12 & 3 	& & 			\\ \hline
456 & 123 & $\emptyset$	& & 			& 123 & $\emptyset$ & & 			& 123 & $\emptyset$ & & 			& 123 & $\emptyset$ & & 			\\  \hline
\end{tabular}
\caption{External and internal activities for a subdivision of $U_{3,6}$} 
\label{tableacti}
}
\end{table}
\end{landscape}

We will divide the proof of Theorem \ref{th:activities} into two lemmas.
\begin{lemma}\label{lemapoly}
Let $B \subseteq \left[n\right]$, $E \subseteq \left[n\right] \setminus B$ and $I \subseteq B$. Let
\begin{eqnarray*}
V(B,E,I) = \{ A \subseteq \left[n\right] : e_A-e_B = e_a-e_b \text{ with } a \in E \text{ and } a>b, & \\
\text{ or with } b \in I \text{ and } a<b \} & 
\end{eqnarray*}
and
\[
P(B,E,I) = \conv \left\{ \frac{e_A + e_B}{2} : A \in V(B,E,I) \right\}.
\]
Then for any matroid $M \in \M$, we have that $Q(M) \cap P(B,E,I) = \emptyset$ if and only if 
\begin{itemize}
\item  $B$ is not a basis of $M$, or
\item $B$ is a basis of $M$ with $E \subseteq E(B)$ and $I \subseteq I(B)$.
\end{itemize}
\end{lemma}

To illustrate this lemma with an example, consider the case $n=4$, $B=\left\{1,3\right\}$, $E=\left\{2\right\}$ and $I=\left\{ 3 \right\}$. Then $V(B,E,I)=\left\{ \left\{1,2\right\}, \left\{2,3\right\} \right\}$. Figure \ref{p(b,e,i)} shows the polytope $P=P(B,E,I)$ inside the hypersimplex, whose vertices are the characteristic vectors of the $2$-subsets of $[4]$. The polytope of the matroid $M_1$ with bases $\mathcal{B}_1= \left\{ \left\{1,2\right\}, \left\{1,4\right\}, \left\{2,3\right\}, \left\{3,4\right\} \right\}$ does not intersect $P$ because $B$ is not a basis of $M_1$. The polytope of the matroid $M_2$ with bases $\mathcal{B}_2= \left\{ \left\{1,3\right\}, \left\{1,4\right\}, \left\{3,4\right\} \right\}$ does not intersect $P$ either, because $B$ is a basis of $M_2$, but $2$ is externally active with respect to $B$ and $3$ is internally active with respect to $B$. Finally, the polytope of the matroid $M_3$ with bases $\mathcal{B}_3= \left\{ \left\{1,3\right\}, \left\{2,3\right\}, \left\{3,4\right\} \right\}$ does intersect $P$, since $B$ is a basis of $M_3$ and $2$ is not externally active with respect to $B$; the intersection point $\frac12(0110+1010)$ ``certifies" this.

\begin{figure}[htp]
\centering
\includegraphics[height=5cm]{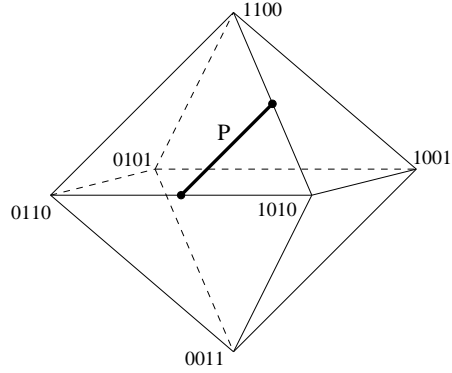}
\caption{The polytope $P=P(B,E,I)$ inside $Q(U_{2,4})$}
\label{p(b,e,i)}
\end{figure}

\begin{proof}
Assume $B$ is a basis of $M$. For $a \notin B$, $a$ is externally active with respect to $B$ if and only if there are no edges in $Q(M)$ of the form $e_a-e_b$ with $a>b$ which are incident to $e_B$. In the same way, for $b \in B$, $b$ is internally active with respect to $B$ if and only if there are no edges in $Q(M)$ of the form $e_a-e_b$ with $a<b$ which are incident to $e_B$. Since the vertices of $P(B,E,I)$ are precisely the midpoints of these edges when $a \in E$ and $b \in I$, if $Q(M) \cap P(B,E,I) = \emptyset$ then $E \subseteq E(B)$ and $I \subseteq I(B)$.

To prove the other direction, suppose that $Q(M) \cap P(B,E,I) \neq \emptyset$. First notice that, since $P(B,E,I)$ is on the hyperplane $x_1+ x_2+ \dotsb + x_n = \left|B\right|$ and $Q(M)$ is on the hyperplane $x_1+ x_2+ \dotsb + x_n = r(M)$, we must have $\left|B\right| = r(M)$. Moreover, since the vertices $v$ of $P(B,E,I)$ satisfy $e_B \cdot v = r(M) - 1/2$ then $B$ must be a basis of $M$, or else the vertices $w$ of $Q(M)$ would all satisfy $e_B \cdot w \leq r(M) - 1$.

Now let $q \in Q(M) \cap P(B,E,I)$. Since $q \in Q(M)$, we know that $q$ is in the cone with vertex $e_B$ generated by the edges of $Q(M)$ incident to $e_B$. In other words, if $A_1, A_2, \dotsc, A_m$ are the bases adjacent to $B$,
\[
q = e_B + \sum_{i=1}^m \lambda_i (e_{A_i} - e_B),  
\]
where the $\lambda_i$ are all nonnegative. If we let $e_{c_i} - e_{d_i} = e_{A_i} - e_B$, then
\[
q = e_B + \sum_{i=1}^m \lambda_i (e_{c_i} - e_{d_i}). 
\]
On the other hand, since $q \in P(B,E,I)$,
\[
q = \sum_{A \in V(B,E,I)} \gamma_A \, \frac{e_A + e_B}{2},
\]
where the $\gamma_A$ are nonnegative and add up to $1$. Setting these two expressions equal to each other we obtain
\[
q = e_B + \sum_{i=1}^m \lambda_i (e_{c_i} - e_{d_i}) = \sum_{A \in V(B,E,I)} \gamma_A \, \frac{e_A + e_B}{2}
\]
and therefore
\[
r = q - e_B = \sum_{i=1}^m \lambda_i (e_{c_i} - e_{d_i}) = \sum_{A \in V(B,E,I)} \gamma_A \, \frac{e_A - e_B}{2}.
\]
For $A \in V(B,E,I)$ we will let $e_{a_A} - e_{b_A}= e_A-e_B$. We have
\begin{equation}\label{sumas}
r = \sum_{i=1}^m \lambda_i (e_{c_i} - e_{d_i}) = \sum_{A \in V(B,E,I)} \gamma_A \, \frac{e_{a_A} - e_{b_A}}{2}.
\end{equation}
Notice that there is no cancellation of terms in either side of  \eqref{sumas}, since the $d_i$s and the $b_A$s are elements of $B$, while the $c_i$s and the $a_A$s are not. Let $r=(r_1, r_2, \dotsc, r_n)$ and let $k$ be the largest integer for which $r_k$ is nonzero.

Assume that $k \notin B$. From the right hand side of \eqref{sumas} and taking into account the definition of $V(B,E,I)$, we have that $k \in E$. From the left hand side we know there is an $i$ such that $c_i=k$. But then $e_{c_i}-e_{d_i}$ is an edge of $Q(M)$ incident to $e_B$, and $d_i < k = c_i$ by our choice of $k$. It follows that $k$ is not externally active with respect to $B$. In the case that $k \in B$, we obtain similarly that $k \in I$, and that $d_j=k$ for some $j$. Thus $e_{c_j}-e_{d_j}$ is an edge of $Q(M)$ incident to $e_B$ and $c_j < k = d_j$, so $k$ is not internally active with respect to $B$. In either case we conclude that $E\nsubseteq E(B)$ or $I \nsubseteq I(B)$, which finishes the proof.
\end{proof}

\begin{lemma}\label{g}
Let $B$ be a subset of $\left[n\right]$, and let $E \subseteq \left[n\right] \setminus B$ and $I \subseteq B$. The function $G_{B,E,I}: \M \to \mathbb{Z}$ defined by
\[
G_{B,E,I}(M)=
\begin{cases}
1 & \text{if } B \text{ is a basis of } M, E=E(B) \text{ and } I=I(B), \\
0 & \text{otherwise}, 
\end{cases}
\]
is a valuation.
\end{lemma}
\begin{proof}
To simplify the notation, we will write $\overline{i_B}$ instead of $\overline{i_{\left\{e_B\right\}}}$. We will prove that $G(B,E,I)=G'(B,E,I)$ where
\begin{equation}\label{g'}
G'_{B,E,I}(M)= (-1)^{\left| E \right| + \left| I \right|} \, \cdot \, \sum_{\substack{E \subseteq X \subseteq \left[n\right] \\ I \subseteq Y \subseteq \left[n\right]}} (-1)^{\left| X \right| + \left| Y \right|} \left( \overline{i_{P(B,X,Y)}}(M) - \overline{i_B}(M) \right),
\end{equation}
which is a sum of valuations.

Let $M \in \M$. If $B$ is not a basis of $M$ then $\overline{i_B}(M)=1$, and by Lemma \ref{lemapoly} we have $\overline{i_{P(B,X,Y)}}(M) = 1$ for all $X$ and $Y$. Therefore $G'_{B,E,I}(M)= 0 =G_{B,E,I}(M)$ as desired. If $B$ is a basis of $M$ then $\overline{i_B}(M)=0$; and we use Lemma \ref{lemapoly} to rewrite \eqref{g'} as
\begin{align*}
G'_{B,E,I}(M) &= (-1)^{\left| E \right| + \left| I \right|} \, \cdot \, \sum_{\substack{E \subseteq X \subseteq E(B) \\ I \subseteq Y \subseteq I(B)}} (-1)^{\left| X \right| + \left| Y \right|} \\
&= (-1)^{\left| E \right| + \left| I \right|} \, \cdot \, \sum_{E \subseteq X \subseteq E(B)} (-1)^{\left| X \right|} \, \cdot \, \sum_{I \subseteq Y \subseteq I(B)} (-1)^{\left| Y \right|} \\
&=
\begin{cases}
1 & \text{if } E=E(B) \text{ and } I=I(B), \\
0 & \text{otherwise}, 
\end{cases}
\end{align*}
as desired.
\end{proof}

\begin{proof}[Proof of Theorem \ref{th:activities}]
The coefficient of $(B,E,I)$ in the definition of \eqref{activity} is $G_{B,E,I}(M)$, so the result follows from Lemma \ref{g}.
\end{proof}

Theorem \ref{th:activities} is significantly stronger than the following result of Speyer which motivated it:

\begin{corollary}(Speyer, \cite{Speyer})
The Tutte polynomial (and therefore any of its evaluations) is a valuation under matroid subdivisions.
\end{corollary}
\begin{proof}
By Theorem \ref{tutteacti}, $T_M(x,y)$ is the composition of the function $h:G \to \mathbb{Z}[x,y]$ defined by $h(B,E,I)=x^{\left| I \right|} \, y^{\left| E \right|}$ with the function $F$ of Theorem \ref{th:activities}.
\end{proof}

\section{\textsf{Related work.}}\label{sec:related work}  

Previous to our work, 
Billera, Jia and Reiner~\cite{BJR}  and Speyer~\cite{Speyer, Speyer2} had studied various valuations of matroid polytopes. A few months after our paper was submitted, we learned about Derksen's results on this topic \cite{Derksen}, which were obtained independently and roughly simultaneously. 
Their approaches differ from ours in the basic fact that 
they are concerned with matroid invariants which are valuations, whereas
our matroid functions are not necessarily constant under matroid isomorphism; however there are similarities. We outline their main invariants here.  

In his work on tropical linear spaces \cite{Speyer}, Speyer shows that the Tutte polynomial is
a valuative invariant.  He also defines in \cite{Speyer2} a polynomial invariant~$g_M(t)$ of a matroid $M$ which arises in the $K$-theory of the Grassmannian.  It is not known how to describe $g_M(t)$ combinatorially in terms of $M$. 

Given a matroid~$M=(E,\B)$, a function $f:E\to\mathbb Z_{>0}$ is
said to be 
{\em $M$-generic} if the minimum value of 
$\sum_{b\in B} f(b)$ over all bases $B\in\B$ is attained just once. 
Billera, Jia, and Reiner study the valuation 
\[
QS(M) = \sum_{\mbox{\scriptsize $f$\ $M$-generic}}
\, \prod_{b\in E} x_{f(b)},
\]
which takes values in the ring of {\em quasi-symmetric functions}
in the variables $x_1, x_2, \ldots$; \emph{i.e.}, the ring generated by
\[\sum_{i_1<\ldots<i_r} x_{i_1}^{\alpha_1}\cdots x_{i_r}^{\alpha_r}\]
for all tuples $(\alpha_1,\ldots,\alpha_r)$ of positive integers.

Derksen's invariant is given by
\[G(M) := \sum_{\mathbf A} U(r_M(A_1)-r_M(A_0), \ldots, 
r_M(A_n)-r_M(A_{n-1}))\] 
where ${\mathbf A}=(A_0,\ldots,A_n)$ ranges over all maximal flags of~$M$,
and 
$$\{U({\bf r}) : \mbox{{\bf r} a finite sequence of nonnegative 
integers}\}$$
is a particular basis for the ring of quasi-symmetric functions.
Derksen's invariant can be defined more generally on polymatroids.
He shows that the Tutte polynomial and
the quasisymmetric function of Billera, Jia and Reiner
are specialisations of~$G(M)$, and asks whether $G(M)$
is universal for valuative invariants in this setting.

For the remainder of this section, $F(M)$ will denote the function
of our Theorem~\ref{th:ranks}. Since $F(M)$ is not a matroid invariant,
it cannot be a specialisation of $g_M(t), QS(M)$, or $G(M)$. 
In the other direction, we suspect that, like the Tutte polynomial,
Speyer's polynomial $g_M(t)$ is a specialisation of $F(M)$. 
As one would expect, $G(M)$ and $QS(M)$ are not
specialisations of $F(M)$.
One linear combination that certifies
this is set out in Table~\ref{tab:syzygy}, in which, to
facilitate carrying out the relevant checks for~$F(M)$, the
relevant matroids are specified via their rank functions.

\begin{table}\small
\begin{center}
\begin{tabular}{|c||c|c|c|c|c|c|c|c|}\hline
$S$ & $\emptyset$ & 1 & 2 & 12 & 3 & 13 & 23 & 123 \\\hline
$r_{M_1}(S)$   & 0&1&1&1&0&1&1&1                 \\\hline 
$r_{M_2}(S)$   & 0&1&0&1&1&1&1&1                 \\\hline 
$r_{M_3}(S)$   & 0&0&1&1&1&1&1&1                 \\\hline 
$r_{M_4}(S)$   & 0&1&1&1&1&1&1&1                 \\\hline 
$r_{M_5}(S)$   & 0&1&1&2&0&1&1&2                 \\\hline 
$r_{M_6}(S)$   & 0&1&0&1&1&2&1&2                 \\\hline 
$r_{M_7}(S)$   & 0&0&1&1&1&1&2&2                 \\\hline 
$r_{M_8}(S)$   & 0&1&1&2&1&2&2&2                 \\\hline 
$r_{M_9}(S)$   & 0&0&1&1&1&1&2&2                 \\\hline 
$r_{M_{10}}(S)$& 0&1&1&2&1&2&2&3                 \\\hline 
$r_{M_{11}}(S)$& 0&1&1&2&1&1&2&2                 \\\hline 
$r_{M_{12}}(S)$& 0&1&1&2&1&2&2&3                 \\\hline 
\end{tabular}

\begin{tabular}{|c||c|c|c|c|c|c|c|c|}\hline
$S$ & 4 & 14 & 24 & 124 & 34 & 134 & 234 & 1234 \\\hline
$r_{M_1}(S)$   &                 0&1&1&1&0&1&1&1 \\\hline
$r_{M_2}(S)$   &                 1&1&1&1&1&1&1&1 \\\hline
$r_{M_3}(S)$   &                 1&1&1&1&1&1&1&1 \\\hline
$r_{M_4}(S)$   &                 1&1&1&1&1&1&1&1 \\\hline
$r_{M_5}(S)$   &                 0&1&1&2&0&1&1&2 \\\hline
$r_{M_6}(S)$   &                 1&2&1&2&1&2&1&2 \\\hline
$r_{M_7}(S)$   &                 1&1&2&2&1&1&2&2 \\\hline
$r_{M_8}(S)$   &                 1&2&2&2&1&2&2&2 \\\hline
$r_{M_9}(S)$   &                 1&1&2&2&2&2&3&3 \\\hline
$r_{M_{10}}(S)$&                 1&1&2&2&2&2&3&3 \\\hline
$r_{M_{11}}(S)$&                 1&2&2&3&2&2&3&3 \\\hline
$r_{M_{12}}(S)$&                 1&2&2&3&2&2&3&3 \\\hline
\end{tabular}

\bigskip

\begin{tabular}{|c||c|c|c|c|c|c|c|c|c|c|c|c|}\hline
$i$   & 1 & 2 & 3 & 4 & 5 & 6 & 7 & 8 & 9 &10 &11 &12 \\\hline
$c_i$ &-1 & 1 &-1 & 1 & 1 &-1 &-1 & 1 & 2 &-2 &-2 & 2 \\\hline
\end{tabular}
\end{center}
\caption{The top table contains the rank functions of twelve matroids $M_i$ on~$[4]$, 
$i=1,\ldots,12$.  The bottom table shows coefficients $c_i$ such that
$\sum c_i F(M_i)=0$ but $\sum c_i G(M_i)\neq0$ and $\sum c_i QS(M_i)\neq0$.}
\label{tab:syzygy}
\end{table}

However, one can give a valuation which is similar in spirit to our $F(M)$ and 
specialises to Derksen's $G(M)$.

\begin{proposition} \label{prop:flags}
The function $H:\M\to G^{n}$ defined by
\[H(M) = \sum_{\mathbf A} \bigl((A_1, r(A_1)), \ldots, (A_n, r(A_n))\bigr),\]
where ${\mathbf A}=(A_1,\ldots,A_n)$ ranges over all maximal flags of~$M$,
is a valuation.
\end{proposition}

\begin{proof}
The proof is a straightforward extension of our argument for
Theorem~\ref{th:ranks}.  With the notation of that proof, 
checking whether a matroid $M$ satisfies $r_M(A_i)=r_i$
for some fixed vector~$\mathbf r=(r_i)$, 
\emph{i.e.}\ whether the term $((A_1,r_1),\ldots,(A_n,r_n))$
is present in $H(M)$, is equivalent to checking that 
$Q(M)$ intersects $P_{A_i,r_i}$ and not $P_{A_i, r_i+1}$ 
for each $i$. 

Observe that if $Q(M)$ intersects $P_{A_i,s_i}$ for all $i$ then
$r(A_i) \geq s_i$ and, since ${\mathbf A}$ is a flag, we can 
choose a single basis of $M$ whose intersection with $A_i$
has at least $s_i$ elements for each $i$. Therefore $Q(M)$ intersects
 $P_{A_1, s_1} \cap \cdots \cap P_{A_n,s_n}$.

Consider the sum
\begin{equation}\label{inexflag}
\sum (-1)^{e_1+...+e_n} i_{P_{\mathbf A, {\mathbf r+\mathbf e}}}(M)
\end{equation}
where the sum is over all $\mathbf e=(e_1, ...,e_n)\in \{0,1\}^{n}$, 
and where $P_{\mathbf A, {\mathbf r+\mathbf e}}$ is the intersection
$P_{A_1, r_1+e_1} \cap \cdots \cap P_{A_n, r_n+e_n}$. 
By our previous observation this sum equals 
\[
\left(\sum_{e_1} (-1)^{e_1} i_{P_{A_1, r_1+e_1}}(M)\right) \dotsb \left(\sum_{e_n} (-1)^{e_n} i_{P_{A_n, r_n+e_n}}(M)\right),
\]
which is $1$ if the term $((A_1,r_1),\ldots,(A_n,r_n))$
is present in $H(M)$, and is $0$ otherwise. 
All the terms in~\eqref{inexflag} are valuations, hence 
$H$ is a valuation.  
\end{proof}


\section{\textsf{Acknowledgments.}}

This work was carried out as part of the SFSU-Colombia Combinatorics Initiative. We are grateful to the San Francisco State University for their financial support of this initiative, and to the Universidad de Los Andes for supporting the third author's visit to SFSU in the Summer of 2007. The third author also wants to thank the Fondo de Investigaciones de la Facultad de Ciencias of the Universidad de los Andes for financially supporting his research.

We would also like to thank the referee for a careful reading of this work; the suggestions in the report led to Table \ref{tab:syzygy} and Proposition \ref{prop:flags}.

\end{document}